\documentclass{gtmon_a}
\pdfoutput=1

\usepackage{amscd, pinlabel}

\proceedingstitle{Groups, homotopy and configuration spaces (Tokyo
  2005)}
\conferencestart{5 July 2005}
\conferenceend{11 July 2005}
\conferencename{Groups, homotopy and configuration spaces, 
                in honour of Fred Cohen's 60th birthday}
\conferencelocation{University of Tokyo, Japan}

\editor{Norio Iwase}
\givenname{Norio}
\surname{Iwase}

\editor{Toshitake Kohno}
\givenname{Toshitake}
\surname{Kohno}

\editor{Ran Levi}
\givenname{Ran}
\surname{Levi}

\editor{Dai Tamaki}
\givenname{Dai}
\surname{Tamaki}

\editor{Jie Wu}
\givenname{Jie}
\surname{Wu}

\title{On braid groups and homotopy groups}

\author{F\,R Cohen}
\givenname{F\,R} 
\surname{Cohen}
\address{Department of Mathematics\\University of Rochester\\\newline
Rochester, NY 14627\\USA} 
\email{cohf@math.rochester.edu}
\urladdr{http://www.math.rochester.edu/people/faculty/cohf/}

\author{J Wu}
\givenname{Jie} \surname{Wu}
\address{Department of Mathematics\\
National University of Singapore\\\newline
Singapore 117543\\
Republic of Singapore} 
\email{matwuj@nus.edu.sg}
\urladdr{http://www.math.nus.edu.sg/~matwujie}


\arxivreference{} 

\volumenumber{13} 
\issuenumber{} 
\publicationyear{2008} 
\papernumber{07}
\startpage{169}
\endpage{193}
\doi{} 
\MR{} 
\Zbl{} 
\received{11 July 2006} 
\revised{29 May 2007} 
\accepted{1 June 2007} 
\published{22 February 2008}
\publishedonline{22 February 2008} 
\proposed{} 
\seconded{} 
\corresponding{}
\version{}

\makeatletter
\def\cnewtheorem#1[#2]#3{\newtheorem{#1}{#3}[section]
\expandafter\let\csname c@#1\endcsname\c@subsection}

\let\xysavmatrix\xymatrix
\def\xymatrix{\disablesubscriptcorrection\xysavmatrix}
\AtBeginDocument{\let\bar\wbar\let\tilde\wtilde\let\hat\what}
\def\<{\langle}
\def\>{\rangle}

\theoremstyle{plain}
\cnewtheorem{prop}[subsection]{Proposition}
\cnewtheorem{thm}[subsection]{Theorem}
\cnewtheorem{lem}[subsection]{Lemma}
\cnewtheorem{cor}[subsection]{Corollary}

\theoremstyle{definition}
\cnewtheorem{defn}[subsection]{Definition}
\cnewtheorem{exm}[subsection]{Example}
\cnewtheorem{remark}[subsection]{Remark}
\makeatother  
\makeautorefname{exm}{Example}

\def\Aut{\mathrm{Aut}}
\def\Brun{\mathrm{Brun}}
\def\QBrun{\mathrm{QBrun}}
\def\Conf{\mathrm{Conf}}
\def\ker{\mathrm{Ker}}
\def\gr{\mathrm{gr}}
\def\Der{\mathrm{Der}}
\def\Sym{\mathrm{Sym}}

\begin{document}

\begin{abstract}    
This article is an exposition of certain connections between the
braid groups, classical homotopy groups of the 2--sphere, as well
as Lie algebras attached to the descending central series of pure
braid groups arising as Vassiliev invariants of pure braids. Natural
related questions are posed at the end of this article.
\end{abstract}

\maketitle


\section{Introduction}\label{sec:Introduction}

The purpose of this article is to give an exposition of certain
connections between the braid groups (Artin \cite{Artin}, see Birman
\cite{Birman}) and classical homotopy groups which arises in joint work
of Jon Berrick, Yan-Loi Wong and the authors
\cite{CW,CW1,BCWW,W}. These connections emerge through several other
natural contexts such as Lie algebras attached to the descending
central series of pure braid groups arising as Vassiliev invariants of
the pure braid groups as developed by T~Kohno \cite{K,K1}.

The main feature of this article is to identify certain
``non-standard'' free subgroups of the braid groups via Vassiliev
invariants of pure braids. A second feature is to indicate natural
ways in which these subjects fit together with classical homotopy
theory. This article is an attempt to draw together these
connections.

Since this paper was submitted, other connections to principal
congruence subgroups in natural matrix groups as well as other
extensions have developed. The authors have taken the liberty of
adding an additional \fullref{sec:Other.connections} with some
of these new connections.

Natural related questions are posed at the end of this article.

Although not yet useful for direct computations, there is a strong
connection between braid groups and homotopy groups. The braids
which naturally arise in this setting also give a large class of
special knots and links arising from Brunnian braids as described
below as well as by Stanford \cite{St}. It is natural to wonder whether and
how these fit together.

The authors take this opportunity to thank Toshitake Kohno, Shigeyuki
Morita, Dai Tamaki as well as other friends for this very enjoyable
opportunity to participate in this conference. The first author is
especially grateful for this mathematical opportunity to learn and to
work on mathematics with friends. The authors also thank an excellent
referee who gave elegant, useful suggestions regarding the exposition
and organization of this article.

\section[Braid groups, Vassiliev invariants and certain free subgroups]{Braid groups, Vassiliev invariants of pure braids and certain free subgroups of braid
groups}\label{sec:Braid groups,Vassiliev invariants of pure braids
and certain free subgroups of braid groups}

The section addresses a naive construction with the braid groups
arising as a ``cabling'' construction. This construction is
interpreted in later sections in terms of the structure of braid
groups, Vassiliev invariants of pure braids as developed by
Toshitake Kohno \cite{K,K1}, associated Lie algebras and the
homotopy groups of the $2$--sphere (Berrick, Wong and the authors 
\cite{CW,CW1,BCWW,W}).

Let $B_k$ denote Artin's $k$--stranded braid group while $P_k$
denotes the pure $k$--stranded braid group, the subgroup of $B_k$
which corresponds to the trivial permutation of the endpoints of the
strands. The group $P_k$ is the fundamental group of the
configuration space of ordered $k$--tuples of distinct points in the
plane $$\Conf(\mathbb R^2,k)$$ for which $\Conf(M,k)= \{(m_1, m_2,
\cdots, m_k)\in M^k\mid  m_i \neq m_j \hbox{ for all}\ i \neq j \}$ for
any space $M$.

The group $B_k$ is the fundamental group of the orbit space
$$\Conf(\mathbb R^2,k)/\Sigma_k$$ obtained from the natural, free
(left)-action of the symmetric group on $k$ letters $\Sigma_k$. The
$k$--stranded braid group of an arbitrary connected surface $S$,
$B_k(S)$, is defined to be the fundamental group of the
configuration space of unordered $k$--tuples of distinct points in
$S$, $\pi_1\Conf(S,k)/\Sigma_k$. The pure braid group $P_k(S)$ is
defined to be the fundamental group $$\pi_1\Conf(S,k).$$
The pure braid groups $P_k$ and $P_k(S^2)$ are intimately related to
the loop space of the 2--sphere as elucidated below in the \fullref{sec:Simplicial objects}. Similar properties are satisfied for
any sphere as described in \fullref{sec:Other.connections}.

To start to address this last point, first consider the free group on
$N$ letters $F_N= F[y_1,\cdots,y_N]$ together with elements $x_i$ for
$ 1 \leq i \leq N $ in $P_{N+1}$ given by the naive ``cabling''
pictured in \fullref{braid} below. The braid $x_1$ with $N=1=i$ in
\fullref{braid} is Artin's generator $A_{1,2}$ of $P_2$. The
braids~$x_i$ for $ 1 \leq i \leq N$ in \fullref{braid} yield
homomorphisms from a free group on $N$ letters $F_N=
F[y_1,\cdots,y_N]$ to $P_{N+1}$
$$\Theta_N\co F[y_1,\cdots,y_N] \to P_{N+1}$$ defined on
generators $y_i$ in $F_N$ by the formula
$$\Theta_N(y_i) = x_i.$$ The maps $\Theta_N$ are the subject of thew
authors' papers \cite{CW,CW1} where it is shown $\Theta_N\co F_N \to
P_{N+1}$ is faithful for every $N$. Three natural questions arise: (1)
Why would one want to know whether $\Theta_n$ is faithful, (2) are
there sensible applications and (3) why is $\Theta_n$ faithful? The
answers to these three questions provide the main content of this
expository article.
\begin{figure}[ht!]
 \centering
 \labellist\small
 \pinlabel $i$ [b] at 62 763
 \pinlabel $N{+}1{-}i$ [b] at 192 763
 \endlabellist
 \includegraphics{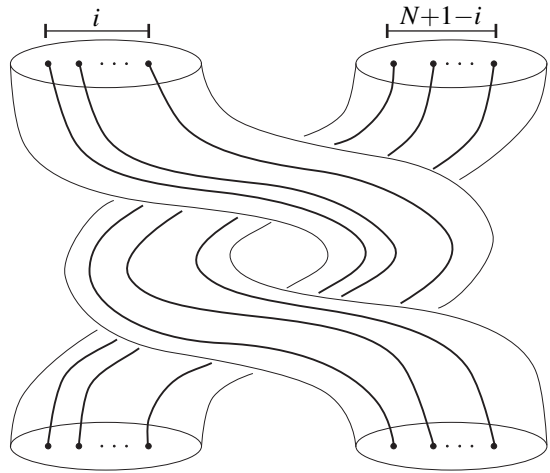}
 \caption{The braid $x_i$ in $P_{N+1}$}
 \label{braid}
\end{figure}

\section{On $\Theta_n$} \label{sec:Why.is.Theta.faithful.}

This section addresses one reason why the map $\Theta_n$ is faithful
\cite{CW,CW1}. The method of proof is to appeal to the structure of the
Lie algebras obtained from the descending central series for both
the source and the target of $\Theta_n$. The structure of these Lie
algebras is reviewed below.

Recall the descending central series of a discrete group $\pi$ given
by $$\pi = \Gamma^1(\pi) \geq \Gamma^1(\pi) \geq \cdots$$ where
$\Gamma^i(\pi)$ is the subgroup of $\pi$ generated by all commutators
$$[\cdots [ x_1,x_2]x_3]\cdots ]x_t]$$ for $ t \geq i$ with $x_i \in
\pi$. The group $\Gamma^i(\pi)$ is a normal subgroup of $\pi$ with the
successive sub-quotients $$gr_i(\pi) =
\Gamma^i(\pi)/\Gamma^{i+1}(\pi)$$ which are abelian groups having
additional structure as follows (cf Magnus, Karrass and Solitar
\cite{MKS}).

Consider the direct sum of all of the $gr_i(\pi)
 = \Gamma^i(\pi)/\Gamma^{i+1}(\pi)$ denoted $$gr_*(\pi) =
\oplus_{i \geq 1}\Gamma^i(\pi)/\Gamma^{i+1}(\pi).$$ The commutator
function
\begin{gather*}[-,-]\co  \pi \times \pi \to \pi\\\tag*{\hbox{given by}}
[x,y] = xyx^{-1}y^{-1},
\end{gather*}
passes to quotients to give a
bilinear map $$[-,-]\co  gr_s(\pi) \otimes_{\Z} gr_t(\pi) \to
gr_{s+t}(\pi)$$ which satisfies both the antisymmetry law and Jacobi
identity for a Lie algebra. (One remark about definitions: The
abelian group $gr_*(\pi)$ is both a graded abelian group and a Lie
algebra, but not a graded Lie algebra as the sign conventions do not
work properly in this context. This situation can be remedied by
doubling all degrees of elements in $gr_*(\pi)$.)

The associated graded Lie algebra obtained from the descending
central series for the target yields Vassiliev invariants of pure
braids by work of Kohno \cite{K,K1}. This Lie algebra has been used
by both Kohno and Drinfel'd \cite{D} in their work on the KZ
equations. The Lie algebra obtained from the descending central
series of the free group $F_N$ is a free Lie algebra by a classical
result due to P~Hall \cite{H}; see also Serre \cite{Serre}.

The proof described next yields more information than just the fact
that $\Theta_N$ is faithful. The method of proof gives a natural
connection of Vassiliev invariants of braids to a classical spectral
sequence abutting to the homotopy groups of the $2$--sphere. Sections
\ref{sec:Pure braid groups of surfaces as Delta groups},
\ref{sec:Brunnian braids} and \ref{sec:Other.connections} below
provide an elucidation of this interconnection.

A discrete group $\pi$ is said to be residually nilpotent provided
$$\bigcap_{i \geq 1}\Gamma^i(\pi) = \{\mathrm{identity}\}$$ where
$\Gamma^i(\pi)$ denotes the $i$-th stage of the descending central
series for $\pi$. Examples of residually nilpotent groups are free
groups, and $P_n$.
\begin{lem}\label{lem:residual.nilpotence}\ 

\begin{enumerate}
  \item Assume that $\pi$ is a residually nilpotent group.
Let $$\alpha\co  \pi  \to\ G$$ be a homomorphism of discrete
groups such that the morphism of associated graded Lie algebras
$$\gr_*(\alpha)\co  \gr_*(\pi) \to\ \gr_*(G)$$ is a monomorphism. Then
$\alpha$ is a monomorphism.
  \item If $\pi$ is a free group, and $\gr_*(\alpha)$ is a
  monomorphism, then $\alpha$ is a monomorphism.
\end{enumerate}
\end{lem}

Thus one step in the proof of \fullref{thm: values of theta}
below is to describe the map $$\Theta_n \co F[y_1,y_2, \cdots,
y_n] \to\ P_{n+1}$$ on the level of associated graded Lie algebras
$$\gr_*(\Theta_n) \co \gr_*(F[y_1,y_2, \cdots,
y_n]) \to\ \gr_*(P_{n+1}).$$ Recall Artin's generators $A_{i,j}$ for
$P_{n+1}$ together with the projections of the $A_{i,j}$ to
$\gr_*(P_{n+1})$ labeled $B_{i,j}$ \cite{CW,CW1}.

\begin{thm}\label{thm: values of theta}
The induced morphism of Lie algebras $$\gr_*(\Theta_n) \co
\gr_*(F[y_1,y_2, \cdots, y_n]) \to\ \gr_*(P_{n+1})$$ satisfies the
formula $$\gr_*(\Theta_n)(y_q) = \Sigma_{1 \leq i \leq n-q+1 < j
\leq n+1} B_{i,j}.$$
\end{thm}

Examples of this theorem are listed next.\eject
\begin{exm}\label{exm: values of theta}\ 

\begin{enumerate}
  \item If $q= 1$, then $$\gr_*(\Theta_n)(y_1) =  B_{1,n+1} + B_{2,n+1} + \cdots +
B_{n,n+1}.$$ Thus if $q= 1$,  and $n = 3$,
$$\gr_*(\Theta_3)(y_1) =  B_{1,4} + B_{2,4} +
B_{3,4}.$$
  \item If $q= 2$, then $$\gr_*(\Theta_n)(y_2) =  (B_{1,n+1} + B_{2,n+1} + \cdots +
B_{n-1,n+1}) + (B_{1,n} + B_{2,n} + \cdots + B_{n-1,n}).$$ Thus if
$q= 2$,  and $n = 3$,
$$\gr_*(\Theta_3)(y_2) =  (B_{1,4} + B_{2,4}) +
(B_{1,3}+ B_{2,3}).$$

\item In general, \begin{align*}
\gr_*(\Theta_n)(y_q)& = V_{n-q+2} + V_{n-q+3} + \cdots +
V_{n+1}\\
\tag*{\hbox{\hspace{28pt}where}}V_{r}&= B_{1,r} + B_{2,r} + \cdots + B_{r-1,r}.
\end{align*}
Thus if $q= 3$,  and $n = 3$,
$$\gr_*(\Theta_3)(y_3) =  B_{1,2} + B_{1,3} + B_{1,4}.$$

\end{enumerate}
\end{exm}

To determine the map of Lie algebras with a more global view, the
structure of the Lie algebra $\gr_*(P_n)$ is useful, and is given as
follows. Let $L[S]$ denote the free Lie algebra over $\mathbb Z$
generated by a set $S$. The next theorem was given in work of Kohno
\cite{K,K1}, and Falk and Randell \cite{FR}.

\begin{thm}\label{thm:Lie.algebras.for.pure.braid.groups}
The Lie algebra $\gr_*(P_n)$ is the quotient of the free Lie algebra
generated by $B_{i,j}$ for $1 \leq i < j \leq n$ modulo the {\it
infinitesimal braid relations} (also called the {\it horizontal 4T
relations} or {\it Yang--Baxter--Lie relations})
$$\gr_*(P_n) = L[B_{i,j}\mid  1 \le i < j \le n]/I$$
where $I$ denotes the $2$--sided (Lie) ideal generated by the
infinitesimal braid relations as listed next:

\begin{enumerate}
\item $[B_{i,j}, B_{s,t}] = 0$, if $\{i, j\}{\cap}\{s, t\} =
\emptyset$
\item $[B_{i,j}, B_{i,s} + B_{s,j}] = 0$
\item $[B_{i,j}, B_{i,t} + B_{j,t}] = 0$
\end{enumerate}
\end{thm}

A computation with these maps gives the following result of
\cite{CW,CW1} for which further connections are elucidated in sections
\ref{sec:Pure braid groups of surfaces as Delta groups} and
\ref{sec:Other.connections}.
\begin{thm}\label{thm: embeddings of Lie algebras}
The maps $\Theta_n \co F[y_1,y_2, \cdots, y_n] \to\ P_{n+1}$ on
the level of associated graded Lie algebras
$$\gr_*(\Theta_n) \co \gr_*(F[y_1,y_2, \cdots, y_n]) \to\ \gr_*(P_{n+1})$$
are monomorphisms. Thus the maps $\Theta_n$ are monomorphisms.
\end{thm}

\begin{remark} Two remarks concerning $\Theta_n$ are given next.

\begin{enumerate}
  \item That $\Theta_n$ is a monomorphism identifies $F[y_1,y_2, \cdots,
y_n]$ as a free subgroup of rank $n$ in $P_{n+1}$. However, there
are other, natural free groups of rank $n$ in $P_{n+1}$.  These
arise from the fibrations of Fadell and Neuwirth  given by
projection maps $p_i\co  \Conf(\R^2,n+1) \to \Conf(\R^2,n)$ which delete
the $i$-th coordinate and have fibre of the homotopy type of an
$n$--fold wedge of circles $\vee_n S^1$ (Fadell and Neuwirth \cite{FN}).

Let $q_i\co  P_{n+1} \to P_n$ denote the map induced by $p_i$ on the
level of fundamental groups. The kernel of $q_i$ is a free group of
rank $n$.

The image of $\Theta_n$ has a contrasting feature: Any composite of
the natural projection maps $d_I\co P_{n+1} \to P_2$ precomposed with
$\Theta_n$,
\[
\begin{CD}
F[y_1,y_2, \cdots, y_n] @>{\Theta_n}>> P_{n+1} @>{d_I}>> P_2,
\end{CD}
\] is a surjection.

  \item The combinatorial behavior of the map $\gr_*(\Theta_n)$ is intricate
even though the definition is elementary as well as natural. For
example, various powers of $2$ arise in the computation of the map
$$\gr_*(\Theta_n)\co \gr_*(F[y_1,y_2, \cdots, y_n]) \to\
\gr_*(P_{n+1})$$ for $n>2$. One example is listed next.
\end{enumerate}
\end{remark}

\begin{exm} $\Theta_3([[[y_1,y_2]y_3]y_2]) = -[[[\gamma_1,\gamma_2]\gamma_3]\gamma_2] +
2[[[\gamma_1,\gamma_3]\gamma_2]\gamma_2] + \delta$ where $\delta$ is
independent of the other terms with $\gamma_1 = B_{1,4} + B_{2,4}+
B_{3,4}$, $\gamma_2 = B_{3,4}$ and $\gamma_3 = B_{2,4} + B_{3,4}$.
At first glance, these elements may appear to be "random``. However,
this formula represents a systematic behavior which arises naturally
from kernels of certain morphisms of Lie algebras.
\end{exm}

The crucial feature which makes the computations effective is the
``infinitesimal braid relations''. In addition, the behavior of the
map $\gr_*(\Theta_n)$ is more regular after restricting to certain
sub-Lie algebras arising in the third stage of the descending
central series \cite{CW,CW1}. Finally, the maps $\Theta_n$ also induce
monomorphisms of restricted Lie algebras on passage to the Lie
algebras obtained from the mod-$p$ descending central series
\cite{CW,CW1}.

\section{Simplicial objects, and $\Delta$--objects}\label{sec:Simplicial objects}

Some basic constructions which are part of an algebraic topologist's
toolkit are described in this section (Moore \cite{Moore}, Curtis
\cite{Curtis}, Bousfield and Kan \cite{BK}, May \cite{May}, and Rourke
and Sanderson \cite{RS}).

One of the great insights in classical homotopy theory, due to Moore
and then `extended' by Kan as well as Rourke and Sanderson is that
not only are homology groups a combinatorial invariant, but so are
homotopy groups. The basic combinatorial framework is that
of a simplicial set developed in \cite{Moore,Curtis,BK,May} and a
$\Delta$--set developed in \cite{RS}  both of which
model the combinatorics of a simplicial complex.

\begin{defn}\label{defn:delta.set}
A {\it $\Delta$--set} is a collection of sets $$K_{\bullet} = \{K_0,
K_1, \cdots\}$$ with functions, {\it  face operations,}
$$d_i\co K_t \to K_{t-1}  \text{ for }    \text{ $0 \leq i \leq t$} $$
which satisfy the identities
$$d_id_j = d_{j-1}d_i$$ if $i < j$. A $\Delta$--object in a small
category $\mathcal C$ is a $\Delta$--set $K_{\bullet} = \{K_0, K_1,
\cdots\}$ with the $K_i$ given by objects in $\mathcal C$ and the
maps $d_i\co K_t \to K_{t-1}$ given by morphisms in $\mathcal C$. Thus,
a {\it $\Delta$--group} is a $\Delta$--set for which all $d_i\co K_t \to
K_{t-1}$ are group homomorphisms.

Let $\Delta[n]$ denote the $n$--simplex with $\delta_i\co  \Delta[n] \to
\Delta[n+1]$ the inclusion of the $i$-th face. Assume that each set
$S_i$ is given the discrete topology unless otherwise stated. The
geometric realization of a $\Delta$--set $K_{\bullet}$ is given by
$$ |K_{\bullet}| = (\amalg K_n \times \Delta[n])/\sim$$
where $\sim$ is the equivalence relation generated by requiring that
if $ x \in K_{n+1} $ and $\alpha \in \Delta[n]$, then
$$(d_i(x), \alpha) \sim  (x, \delta_i(\alpha)).$$

\end{defn}

\begin{exm}\label{exm:delta.set}

Examples of $\Delta$--sets are given next.

\begin{enumerate}
   \item A choice of $\Delta$--set with exactly one 0--simplex and one 1--simplex 
has geometric realization given by the circle.

\item  A natural example of a $\Delta$--group arises from the
pure braid groups $P_n(S) = \pi_1\Conf(S,n)$ for a path-connected
surface $S$ \cite{BCWW}. Define
$$\Delta_{\bullet}(S)\quad\text{by}\quad\Delta_n(S) = P_{n+1}(S)$$ the $(n+1)$-st pure 
braid group for the surface $S$.

There are $n+1$ homomorphisms
$$d_i\co P_{n+1}(S) \to P_{n}(S)$$ for $ 0 \leq i \leq n$ obtained by
deleting the $(i+1)$-st strand in $P_{n+1}(S)$. The homomorphisms
$d_i$ are induced on the level of fundamental groups of
configuration spaces by the projection maps
$$p_{i+1}\co \Conf(S, n+1) \to \Conf(S, n)$$ given be deleting the
$(i+1)$-st coordinate and satisfy the identities required for a
$\Delta$--group.

In case $S =\mathbb C \mathbb P^1 = S^2$, the associated
$\Delta$--group gives basic information about the homotopy groups of
the 2--sphere \cite{BCWW}. In case $S_g$ is a closed orientated
surface, the $\Delta$--group $\Delta_{\bullet}(S_g)$ does not admit
the structure of a simplicial group as given in the next \fullref{defn:simplicial.set}.
\end{enumerate}

\end{exm}

Simplicial sets are defined next.
\begin{defn}\label{defn:simplicial.set} A {\it simplicial set}
is

\begin{enumerate}
  \item a $\Delta$--set $K_{\bullet} = \{K_0, K_1, \cdots\}$ together with
  \item functions, {\it  degeneracy operations,} $$s_j\co K_t \to K_{t+1} \text{ for }  \text{ $0 \leq j \leq t$} $$
  which satisfy the {\it simplicial identities }
$$ d_id_j = d_{j-1}d_i  \text{ if $i<j$,}$$
$$s_is_{j} = s_{j+1}s_i  \text{ if $i \leq j$ and}$$
\[
d_is_j =
\begin{cases}
s_{j-1}d_i & \text{if $i<j$,}\\
\mathrm{identity} & \text{if $i=j$ or $i=j+1$,}\\
s_jd_{i-1} & \text{if $i>j+1$.}
\end{cases}
\]
\end{enumerate}

A simplicial-object in a small category $\mathcal C$ is a
simplicial-set $K_{\bullet} = \{K_0, K_1, \cdots\}$ with the $K_i$
given by objects in $\mathcal C$ for which both face maps $d_i\co K_t
\to K_{t-1}$ and degeneracies $s_j\co K_t \to K_{t+1}$ are given by
morphisms in $\mathcal C$.

Thus, a {\it simplicial-group} $$G_{\bullet} = \{G_0,G_1, \cdots
\}$$ is a simplicial-set for which all of the $G_i$ are groups with
face and degeneracies given by group homomorphisms.
\end{defn}

\begin{exm}\label{exm:basic.simplicial.set}
Two examples of simplicial sets are given next.
\begin{enumerate}

\item The simplicial 1--simplex $\Delta[1]$ has two 0--simplices $\<0\>$ and $\<1\>$. The $n$--simplices
of $\Delta[1]$ are sequences $\<0^i, 1^{n+1-i}\>$ for $0 \leq i \leq
n+1$. All of the non-degenerate simplices are $\<0\>$, $\<1\>$, and
$\<0,1\>$.

\item The simplicial circle $S^1$ is a quotient of the simplicial 1--simplex $\Delta[1]$
obtained by identifying $\<0\>$ and $\<1\>$. There are exactly two
equivalence classes of non-degenerate simplices given by $\<0\>$, and
$\<0,1\>$. Furthermore, the simplicial circle $S^1$ is given in degree
$k$ by

\begin{enumerate}
  \item a single point $\<0\>$ in case $k = 0$, and
  \item $n+1$ points $\<0^i, 1^{n+1-i}\>$ for $0 \leq i < n+1$ in case $k = n$
  for which $\<0^{n+1}\>$ and $\<1^{n+1}\>$ are identified.
\end{enumerate}

In what follows below, it is useful to label these simplices by
$$y_{n+1-i} = \<0^i, 1^{n+1-i}\>$$ for $ 0 <i \leq n+1$ with $y_0 =
s_0^n(\<0\>)$.

In addition, classical, elegant constructions for the standard
simplicial $n$--simplex $\Delta[n]$ as well as the $n$--sphere are
given in \cite{BK,Curtis,May}.
\end{enumerate}

\end{exm}

Homotopy groups are defined for any $\Delta$--set, but the definition
admits a simple description for simplicial sets which satisfy an
additional condition.

\begin{defn}\label{defn:extension.condition}
A simplicial set $K_{\bullet}$ is said to satisfy the {\it extension
condition} if for every set of $n$--simplices $x_0$, $x_1$, \dots,
$x_{k-1}$, $x_{k+1}$, \dots, $x_{n+1} \in K_{n}$ which satisfy the
compatibility condition
\begin{equation*}
d_i x_j = d_{j-1} x_i, \ \ i<j, \ \ i \not= k, \ \  j\not= k,
\end{equation*}
there exists an $(n+1)$--simplex $x \in K_{n+1}$ such that $d_i x =
x_i$ for $i \not=k$. A simplicial set which satisfies the extension
condition is also called a {\it Kan complex}.

For each Kan complex $K_{\bullet}$ with base-point $* \in K_0$ and
$s_0^n(*) \in K_n$, the $n$-th homotopy set $\pi_n(K_{\bullet},*)$ is
defined for $n \geq 1$ as the equivalences classes of simplices $a
\in K_n$, denoted $[a]$, such that

\begin{enumerate}
   \item $d_i(a) = s_0^{n-1}(*)$ for all $ 0 \leq i \leq n$, and
   \item two such simplices $a,a' \in K_n$ are equivalent provided
   there exists a simplex $y \in K_{n+1}$ such that $d_{n+1}y = a'$,
$d_{n}y = a$, and $d_{i}y = s_0^n(*)$ for all $ 0 \leq i <n$.
\end{enumerate}
\end{defn}

If $K_{\bullet}$ satisfies the extension condition, and $n\geq 1$,
then $\pi_n(K_{\bullet},*)$ is a group. In case $n\geq 2$,
$\pi_n(K_{\bullet},*)$ is also abelian.

Basic examples of (i) Kan complexes as well as (ii) $\Delta$--groups
which do not admit the structure of a Kan complex are given next.
\begin{exm}\label{exm:extension.condition}
A simplicial group $G_{\bullet} = \{G_0,G_1, \cdots \}$ always
satisfies the extension condition as shown in \cite{Moore}.
\begin{enumerate}
  \item If $G_{\bullet} = \{G_0,G_1, \cdots \}$ is a simplicial group, then the
$n$-th homotopy group of $G_{\bullet}$ is the quotient group
\begin{gather*}
\pi_n(G_{\bullet}) = Z_n/B_n\\\tag*{\hbox{\hspace{28pt}for which}} 
Z_n = \bigcap_{0 \leq i\leq n}\ker( d_i\co G_n \to G_{n-1}),\\
\tag*{\hbox{\hspace{28pt}and}}B_n = d_0\left(\bigcap_{1
\leq i\leq k+1} \ker( d_i\co G_{n+1} \to G_{n})\right).\end{gather*}

  \item Additional, related information is stated next for certain
$\Delta$--groups $\Delta_{\bullet} = \{\Delta_0,\Delta_1, \cdots \}$
which are not necessarily simplicial groups.

    \begin{enumerate}
    \item In this case, the set of left cosets $\pi_n(\Delta_{\bullet} ) = Z_n/B_n$ is still well-defined.
    \item The left cosets $Z_n/B_n$ may not admit a
    natural structure as a group. One case occurs for $S_g$ a closed, oriented surface with
    $\Delta_{\bullet}  = \Delta_{\bullet} (S_g)$ as given in \fullref{exm:delta.set} and in \cite{BCWW}.
    \item The special case of $S = S^2$ has further properties given
    in in \cite{BCWW} and below.
    \end{enumerate}

\end{enumerate}
\end{exm}

An example of a simplicial group obtained naturally from Artin's
pure braid groups is described next.
\begin{exm}\label{exm:AP}

Consider $\Delta$--groups with $\Delta_n(S) = P_{n+1}(S)$ as given in
\fullref{exm:delta.set} for surfaces $S$. Specialize to the
surface $$S = \mathbb R^2.$$
In this case, there are also $n+1$ homomorphisms
$$s_i\co P_{n+1} \to P_{n+2}$$ obtained by ``doubling'' the
$(i+1)$-st strand. The homomorphisms $s_i$ are induced on the level
of fundamental groups by the maps for configuration spaces
$$\mathbb S_{i}\co \Conf(\mathbb R^2, n+1) \to \Conf(\mathbb R^2,
n+2 )$$ defined by the formula $$\mathbb S_{i}(x_1, \cdots, x_{n+1})
= (x_1, \cdots, x_{i+1}, \lambda(x_{i+1}), x_{i+2},\cdots,
x_{n+1})$$ where $\lambda(x_{i+1}) = x_{i+1} + (\epsilon,0)$ for
$(\epsilon,0)$ a point in $\mathbb R^2$ with $$\epsilon = (1/2)\cdot
\min_{t \neq i+1}||x_{i+1}-x_t||.$$ The homomorphisms $d_i$ and
$s_j$ satisfy the simplicial identities \cite{CW,CW1,BCWW}.

Thus the pure groups in case $S = \R^2$ provide an example of a
simplicial group denoted $$\mathrm{AP}_{\bullet}\quad\hbox{with}
\quad\mathrm{AP}_n = P_{n+1}$$ for $n = 0,1,2,3, \cdots $.

\end{exm}

Consider a pointed topological space $(X,*)$. The pointed loop space
of $X$, $\Omega(X)$, has a natural multiplication coming from "loop
sum`` which is not associative, but homotopy associative. Milnor
proved that the loop space of a connected simplicial complex is
homotopy equivalent to a topological group \cite{Milnor2}. James
\cite{James} proved that the loop space of the suspension of a
connected CW--complex is naturally homotopy equivalent to a free
monoid as explained by Hatcher \cite[page 282]{Hatcher}.  Milnor
realized that the James construction could be translated directly into
the the language of simplicial groups as described next \cite{Milnor}.

\begin{defn}\label{defn:F[K]}
Let $K_{\bullet}$ denote a pointed simplicial set (with base-point
$* \in K_0$ and $s_0^n(*) \in K_n$). Define {\it Milnor's free
simplicial group $F[K]_{\bullet}$ } in degree $n$ by
$$F[K]_{n} = F[K_n]/s_0^n(*) = 1.$$ Then $F[K]_{\bullet}$ is a simplicial group with face and degeneracy
operations given by the natural multiplicative extension of those
for $K_{\bullet}$. In addition, the face and degeneracy operations
applied to a generator is either another generator or the identity
element.
\end{defn}

\begin{exm}\label{exm:F[K]}

An example of $F[K]_{\bullet}$ is given by $K_{\bullet} =
S^1_{\bullet}$ the simplicial circle. Notice that $F[S^1]_{\bullet}$
in degree $n$ is isomorphic to the free group $F[S^1]_{n} =
F[y_1,\cdots,y_n]$ by \fullref{exm:basic.simplicial.set}.

\end{exm}

Milnor defined the geometric realization of a simplicial set
$K_{\bullet} = \{K_0, K_1, \cdots\}$ for which the underlying
topology of $K_{\bullet}$ is discrete \cite{Milnor3}. Recall the
inclusion of the $i$-th face $\delta_i\co\Delta[n-1] \to
\Delta[n] $ togther with the projection maps to the $j$-th face
$\sigma_j\co\Delta[n+1] \to \Delta[n]$ \cite{BK,Curtis,May}.

\begin{defn}\label{defn:geometric.realization.of.a.simplicial.set}
The {\it geometric realization of $K_{\bullet}$} is
$$|K_{\bullet}| = (\amalg K_n \times \Delta[n])/\sim$$ where $\sim$ denotes the equivalence
relation generated by requiring
\begin{enumerate}
\item if $ x \in K_{n+1} $ and $\alpha \in \Delta[n]$, then
$(d_i(x), \alpha) \sim  (x, \delta_i(\alpha)),$ and

\item if $ y \in K_{n} $ and $\beta \in \Delta[n+1]$, then
$(x, \sigma_j(\beta)) \sim  (s_j(x), \beta).$
\end{enumerate}

\end{defn}
\begin{thm}\label{thm:Milnor.free.group}
If $K_{\bullet}$ is a reduced simplicial set (that is $K_0$ is
equal to a single point  $\{*\}$), then the geometric realization
$|F[K]_{\bullet}|$ is homotopy equivalent to $\Omega \Sigma
|K_{\bullet}|$. Thus the homotopy groups of $F[K]_{\bullet}$ (as
given in \cite{Moore} and \fullref{exm:extension.condition})
are isomorphic to the homotopy groups of the space $\Omega \Sigma
|K_{\bullet}|$.
\end{thm}

\begin{exm}\label{exm:F[circle]}
Consider the special case of $K_{\bullet} =S^1_{\bullet}$. Then the
geometric realization $|F[S^1]_{\bullet}|$ is homotopy equivalent to
$\Omega S^2$, and there are isomorphisms
$$\pi_nF[S^1]_{\bullet} \to \pi_n\Omega S^2 \cong \pi_{n+1}S^2 .$$

\end{exm}

A partial synthesis of this information is given sections
\ref{sec:Pure braid groups of surfaces as Delta groups},
\ref{sec:Brunnian braids} and \ref{sec:Other.connections}.

\section{Pure braid groups of surfaces as simplicial groups and $\Delta $--groups}\label{sec:Pure braid groups of surfaces as Delta groups}

The homomorphism $\Theta_n \co F[y_1,y_2, \cdots, y_n] \to\
P_{n+1}$ which arises from the cabling operation described in
\fullref{braid} satisfies the following  properties.

\begin{enumerate}
  \item The homomorphisms $\Theta_n \co F[y_1,y_2,
\cdots, y_n] \to\ P_{n+1}$ give a morphism of simplicial groups
$$\Theta\co  F[S^1]_{\bullet} \to\ \mathrm{AP}_{\bullet}$$ for which the homomorphism
$\Theta_n$ is the restriction of $\Theta$ to $F[S^1]_n$.

  \item By \fullref{thm: embeddings of Lie algebras}, the homomorphisms $\Theta_n \co F[y_1,y_2, \cdots, y_n] \to\
P_{n+1}$ are monomorphisms and so the morphism $\Theta\co  F[S^1]
\to\ \mathrm{AP}_{\bullet}$ is a monomorphism of simplicial groups.
  \item There is exactly one morphism of simplicial groups $\Theta$
  with the property that $\Theta_1(y_1) =  A_{1,2}$.
\end{enumerate}

Thus, the picture given in \fullref{braid} is a description for
generators of $F[S^1]_{n}$ in the simplicial group
$F[S^1]_{\bullet}$. These features are summarized next.

\begin{thm}\label{thm:Embedding braids}

The homomorphisms $\Theta_n \co F[y_1,y_2, \cdots, y_n] \to\
P_{n+1}$ (``pictured'' in \fullref{braid}) give the unique
morphism of simplicial groups
$$\Theta\co  F[S^1]_{\bullet} \to\ \mathrm{AP}_{\bullet}$$ with $\Theta_1(y_1) =  A_{1,2}$.
The map $\Theta$ is an embedding. Hence the $n$-th homotopy group of
$F[S^1]$, isomorphic to $\pi_{n+1}(S^2)$, is a natural sub-quotient
of $\mathrm{AP}_{\bullet}$. Furthermore, the smallest sub-simplicial
group of $\mathrm{AP}_{\bullet}$ which contains the element
$\Theta_1(y_1) = A_{1,2}$ is isomorphic to $F[S^1]_{\bullet}$.
\end{thm}

On the other-hand, the homotopy sets for  the $\Delta$--group
$\Delta_{\bullet}(S^2)$ are also giving the homotopy groups of the
2--sphere via a different occurrence of $F[S^1]_{\bullet}$. The
homeomorphism of spaces $$\Conf(S^2,k) \to PGL(2,\C) \times \Conf(S^2
- Q_3,k-3)$$ for $ k \geq 3 $ and where $Q_3$ denotes a set of there
distinct points in $S^2$ is basic for the next Theorem \cite{BCWW}.
\begin{thm} \label{thm:two.sphere}
If $S = S^2$ and $n\geq 4$, then there are isomorphisms
$$\pi_n(\Delta_{\bullet}(S^2)) \to \pi_n(S^2).$$
\end{thm}

The descriptions of homotopy groups implied by these Theorems admit
interpretations in terms of classical, well-studied features of the
braid groups as given in the next section. An extension to all
spheres is given in \cite{CW,CW1} as pointed out in \fullref{sec:Other.connections}.

\section{Brunnian braids, ``almost Brunnian'' braids, and homotopy groups}\label{sec:Brunnian braids}

The homotopy groups of a simplicial group, or the homotopy sets of a
$\Delta$--group admit a combinatorial description as pointed out in
\fullref{exm:extension.condition}. These homotopy sets are the
set of left cosets $Z_n/B_n$ where $Z_n$ is the group of $n$--cycles
and $B_n$ is the group of $n$--boundaries for the $\Delta$--group.

Recall \fullref{exm:delta.set} in which the $\Delta$--group
$\Delta_{\bullet}(S)$ is specified by $\Delta_n(S) = P_{n+1}(S)$ the
$(n+1)$--stranded pure braid group for a connected surface $S$. The
main point of this section is that the $n$--cycles $Z_n$ are given by
the "Brunnian braids`` while the $n$--boundaries $B_n$ are given by
the ``almost Brunnian braids'', subgroups considered next which are
also important in other applications (Magnum and Stanford
\cite{Mangum.Stanford}).
\begin{defn}\label{defn:brunnian.almostbrunnian}

Consider the $n$--stranded pure braid group for any (connected)
surface $S$, the fundamental group of $\Conf(S,n)$. The group of
Brunnian braids $\Brun_n(S)$ is the subgroup of $P_n(S)$ given by
those braids which become trivial after deleting any single strand.
That is, $$\Brun_n(S) = \bigcap_{0 \leq i\leq n-1} \ker( d_i\co
P_n(S) \to P_{n-1}(S)).$$
The ``almost Brunnian'' $(n+1)$--stranded braid group is
$$\QBrun_{n+1}(S) = \bigcap_{1 \leq i\leq n} \ker( d_i\co
P_{n+1}(S) \to P_{n}(S)).$$ The subgroup $\QBrun_{n+1}(S)$ of
$P_{n+1}(S)$ consists of those braids which are trivial after
deleting any one of the strands $2, 3, \cdots, n+1$, but not
necessarily the first.
\end{defn}

\begin{exm}\label{exm:brunnian.almost.brunnian.AP}
Consider the simplicial group $\mathrm{AP}_{\bullet}$ with
$$\mathrm{AP}_n = P_{n+1}$$ for $n = 0,1,2,3, \cdots $ as given in
\fullref{exm:AP}.

In this case, notice that that the map $d_0\co \QBrun_{k+2} \to
\Brun_{k+1}$ is a split epimorphism. Thus the homotopy groups of the
simplicial group $\mathrm{AP}_{\bullet}$ are all trivial.
\end{exm}

An inspection of definitions gives that
\begin{enumerate}
  \item the group of $n$--cycles of $\Delta_{\bullet}(S)$, $Z_n(S)$, is
precisely $\Brun_{n+1}(S)$ while
  \item the group of $n$--boundaries,
$B_n(S)$, is exactly $d_0(\QBrun_{n+2}(S))$.
\end{enumerate} This feature is recorded next as a lemma.
\begin{lem}\label{lem:cycles.boundaries.Brunnian.braids}
Let $S$ denote a connected surface with associated $\Delta$--group
$\Delta_{\bullet}(S)$ (as given in \fullref{exm:delta.set}).
Then the following hold.

\begin{enumerate}
\item The group of $n$--cycles $Z_n(S)$ is $\Brun_{n+1}(S)$.
\item The group of boundaries $B_n(S)$ is $d_0(\QBrun_{n+2}(S))$.
\item There is an isomorphism
$$\pi_k(\mathrm{AP}_{\bullet}) \to \Brun_{k+1}/d_0(\QBrun_{k+2}).$$
Furthermore, $\pi_k(\mathrm{AP}_{\bullet})$ is the trivial group.

\item There is an isomorphism of left cosets which is natural
for pointed embeddings of connected surfaces $S$
$$\pi_k(\Delta_*(S)) \to \Brun_{k+1}(S)/d_0(\QBrun_{k+2}(S)).$$
\end{enumerate}
\end{lem}

Properties of the $\Delta$--group for the $2$--sphere $S = \mathbb C
\mathbb P^1 = S^2$ is the main subject of \cite{BCWW} where the next
result is proven.
\begin{thm} \label{thm:two.sphere.returned}
If $S = S^2$ and $k\geq 4$, then
$$\pi_k(\Delta_{\bullet}(S^2)) = \Brun_{k+1}(S^2)/d_0(\QBrun_{k+2}(S^2))$$ is a group which is
isomorphic to the classical homotopy group $\pi_k(S^2)$.

Furthermore, there is an exact sequence of groups
$$1 \to \Brun_{k+2}(S^2) \to \Brun_{k+1}(\mathbb R^2) \to  \Brun_{k+1}(S^2)
\to \pi_{k}(S^2) \to 1.$$
\end{thm}

The next lemma follows by a direct check of the long exact homotopy
sequence obtained from the Fadell--Neuwirth fibrations for
configuration spaces \cite{FN,FH}.
\begin{lem}\label{lem:Brunnian.braids}
If $S$ is a surface not homeomorphic to either $S^2$ or $\R\mathbb
P^2$, and $k \geq 3$, then $\Brun_k(S)$ and $\QBrun_{k}(S)$ are free
groups. If $S$ is any surface, and $k \geq 4$, then $\Brun_k(S)$ and
$\QBrun_{k}(S)$ are free groups.
\end{lem}

\begin{exm}\label{exm:example.of.brunnian.braid.group}
One classical example of a Brunnian braid group is $\Brun_4(S^2)$
which is isomorphic to the principle congruence subgroup of level
$4$ in $PSL(2,\Z)$ as given in \fullref{sec:Other.connections}
below.
\end{exm}

One question below in \fullref{sec:Questions} is to consider
properties of the free groups obtained from the intersections
$\Theta_k(F_{k})\cap \Brun_{k+1}$ as well as $\Theta_k(F_{k})\cap
d_0(\QBrun_{k+2})$ where $\Theta_k \co F[y_1,y_2, \cdots, y_k]
\to\ P_{k+1}$ is the homomorphism in \fullref{sec:Why.is.Theta.faithful.}. These groups are precisely the
cycles and boundaries for $F[S^1]_{\bullet}$.

\begin{lem}\label{lem:free.intersect.Brunnian.braids}
If $k \geq 3$, then $\Theta_k(F_{k})\cap \Brun_{k+1}$ as well as
$\Theta_k(F_{k})\cap d_0(\QBrun_{k+2})$ are countably infinitely
generated free groups.
\end{lem}

The standard Hall collection process or natural variations can be used
to give inductive recipes rather than closed forms for
generators. T~Stanford has given a related elegant exposition of the
Hall collection process \cite{St}. The analogous process was applied
by Cohen and Levi~\cite{cohen.levi} to give group theoretic models for
iterated loop spaces.

The connection of the homotopy groups of $S^2$ as well as the Lie
algebra attached to the descending central series of the pure braid
groups is discussed next.

\begin{thm}\label{thm:homotopy.and.embedding braids}
The group $\Theta_k(F_{k})\cap d_0(\QBrun_{k+2})$ is a normal
subgroup of $\Theta_k(F_{k})\cap \Brun_{k+1}$. There are
isomorphisms $$\Theta_k(F_{k})\cap \Brun_{k+1}/\Theta_k(F_{k})\cap
d_0(\QBrun_{k+2}) \to \pi_{k+1}S^2.$$
\end{thm}

The method of proving that the maps $\Theta_n \co F[y_1,y_2,
\cdots, y_n] \to\ P_{n+1}$ are monomorphisms via Lie algebras admits
an interpretation in terms of classical homotopy theory. The method
is to filter both simplicial groups $F[S^1]_{\bullet}$ and
$\mathrm{AP}_{\bullet}$ via the descending central series, and then
to analyze the natural map on the level of associated graded Lie
algebras.

On the other-hand, the Lie algebra arising from filtering any
simplicial group by its' descending central series gives the
$E^0$--term of the Bousfield--Kan spectral sequence for the simplicial
group in question \cite{BK}. Similarly, filtering via the mod-$p$
descending central series gives the classical unstable Adams
spectral sequence \cite{BK,Curtis,W}.

Thus the method of proof of \fullref{thm: embeddings of Lie
algebras} is precisely an analysis of the behavior of the natural
map  $\Theta\co  F[S^1]_{\bullet} \to \mathrm{AP}_{\bullet}$ on the
level of the $E^0$--term of the Bousfield--Kan spectral sequence. This
method exhibits a close connection between Vassiliev invariants of
pure braids and these natural spectral sequences. The next result is
restatement of \fullref{thm: embeddings of Lie algebras} proven
in \cite{CW,CW1}.
\begin{cor}\label{cor: unext.and.embeddings of Lie algebras}
The maps $\Theta_n \co F[y_1,y_2, \cdots, y_n] \to\ P_{n+1}$ on
the level of associated graded Lie algebras
$$\gr_*(\Theta_n) \co \gr_*(F[y_1,y_2, \cdots, y_n]) \to\ \gr_*(P_{n+1})$$
are monomorphisms.  Thus the maps $\Theta_n$ induce embeddings on
the level of the $E^0$--term of the Bousfield--Kan spectral sequences
for $E^0(\Theta) \co E^0(F[S^1]_{\bullet})\to
E^0(\mathrm{AP}_{\bullet})$.
\end{cor}

\section{Other connections}\label{sec:Other.connections}

Several further connections, outlined next, emerged after this
article was submitted.

\medskip
{\bf Connection to principal congruence subgroups}

One basic construction above is the Brunnian braid groups
$\Brun_k(S)$.  Recently, the authors have proven (unpublished) that
the Brunnian braid group $\Brun_4(S^2)$ is isomorphic to the principal
congruence subgroup of level $4$ in $PSL(2,\Z)$ \cite{BCWW}.

This identification may admit an extension by considering the
Brunnian braid groups $\Brun_{2g}(S^2)$ as natural subgroups of
mapping class groups for genus $g$ surfaces. The subgroups
$\Brun_{2g}(S^2)$ may embed naturally in $Sp(2g, \Z)$ via classical
surface topology using branched covers of the 2--sphere (work in
progress).

\medskip
{\bf Connections to other spheres}

The work above has been extended to all spheres as well as other
connected CW--complexes \cite{CW1}. One way in which other spheres
arise is via the induced embedding of free products of simplicial
groups $$ \Theta \amalg \Theta\co F[S^1]_{\bullet} \amalg
F[S^1]_{\bullet} \to \mathrm{AP}_{\bullet}\amalg
\mathrm{AP}_{\bullet}.$$ The geometric realization of
$F[S^1]_{\bullet} \amalg F[S^1]_{\bullet}$ is homotopy equivalent to
$\Omega(S^2 \vee S^2)$ by Milnor's theorem stated above as
\ref{thm:Milnor.free.group}.

Furthermore, $\Omega(S^2 \vee S^2)$ is homotopy equivalent to a weak
infinite product of spaces $\Omega(S^n)$ for all $n > 1$.

\medskip
{\bf Connection to certain Galois groups}

Consider automorphism groups $\Aut(H)$ where $H$ is $F_n$ or a
completion of $F_n$ given by either the pro-finite completion
$\widehat{F_n}$ or the pro-$\ell$ completion
${\widehat{(F_n)}}_{\ell}$.  Certain Galois groups $G$ are identified
as natural subgroups of these automorphism groups by Bely{\u\i}
\cite{Belyi}, Deligne \cite{Deligne}, Drinfel'd \cite{D,D2},
Ihara \cite{Ihara,Ihara2} and Schneps \cite{Schneps}.

One example is Drinfel'd's Grothendieck--Teichm\"uller Galois group
$G = \widehat{GT}$, a subgroup of $\Aut(\widehat{F_2})$.

Let $\Der(L^{R}[V_n])$ denote the Lie algebra of derivations of the
free Lie algebra $L^{R}[V_n]$ where $V_n$ denotes a free module of
rank $n$ over $R$ a commutative ring with identity. Two natural
morphisms of Lie algebras which take values in $\Der(L^{R}[V_n])$
occur in this context as follows.

Properties of the infinitesimal braid relations as stated in \fullref{thm:Lie.algebras.for.pure.braid.groups} above give a natural
second map $$Ad\co gr_*(P_{n+1})\to \Der(L^{\mathbb Z}[V_n])$$ for
which the kernel of $Ad$ is precisely the center of $gr_*(P_{n+1})$
(Cohen and Prassidis \cite{cohen.prassidis}). Combining this last fact
with \fullref{thm: values of theta} gives properties of the
composite morphism of Lie algebras
\[
\begin{CD}
gr_*(F_n)  @>{\gr_*(\Theta_n)}>> gr_*(P_{n+1}) @>{Ad}>>
\Der(L^{\mathbb Z}[V_n]).
\end{CD}
\]

\begin{prop}\label{prop:another.embedding.DerL}
If $ n \geq 2$, the induced morphism of Lie algebras
$$Ad\circ\gr_*(\Theta_n) \co \gr_*(F[y_1,y_2, \cdots, y_n]) \to\
\Der(L^{\mathbb Z}[V_n])$$ is a monomorphism.
\end{prop}

In addition, certain Galois groups $G$ above are filtered with
induced morphisms of Lie algebras $$gr_*(G) \to
\Der(L^{\widehat{\mathbb Z}}[V_n])$$ where $\widehat{\mathbb Z}$
denotes the pro-finite completion of the integers. One example is $G
= \widehat{GT}$ with
$$gr_*(\widehat{GT}) \to \Der(L^{\widehat{\mathbb Z}}[V_2])$$
as given in in \cite{Deligne,Ihara,Ihara2,Schneps2}.

This raises the question of (i) whether the images of
$Ad\circ\gr_*(\Theta_2)$ and $gr_*(\widehat{GT})$ in
$\Der(L^{\widehat{\mathbb Z}}[V_2])$ have a non-trivial intersection
and (ii) whether these intersections are meaningful.

\section{Questions}\label{sec:Questions}

The point of this section is to consider whether the connections
between the braid groups and homotopy groups above are useful. Some
natural as well as speculative problems are listed next.
\begin{enumerate}
\item The combinatorial problem of distinguishing elements in the
pure braid groups has been well-studied. For example, the Lie
algebra associated to the descending central series of the pure
braid group $P_n$ has been connected with Vassiliev theory and has
been shown to be a complete set of invariants which distinguish all
elements in $P_n$ \cite{K1}. Furthermore, these Lie algebras have
been applied to other questions arising from the classical
KZ--equations \cite{K,D} as well as the structure of certain Galois
groups \cite{Ihara,Deligne,D,D2}.

Find invariants of braids up to the coarser equivalence relation
given in \fullref{thm:homotopy.and.embedding braids}. In
particular, can one identify the subset of Vassiliev invariants of
braids which are "homotopy invariant"?

\item
 Give combinatorial properties of the natural map
$\Brun_{k+1}(\mathbb R^2) \to  \Brun_{k+1}( S^2)$  which provide
information about the cokernel. Two concrete problems are stated
next.

\begin{enumerate}
  \item Give group theoretic reasons why the order of the $2$--torsion in
$\pi_*(S^2)$ is bounded above by $4$ and why the $p$--torsion for an
odd prime $p$ is bounded above by $p$.
  \item If $k+1 \geq 5$, the image of $\Brun_{k+1}(\mathbb R^2) \to
\Brun_{k+1}( S^2)$ is a normal subgroup of finite index.

This fact follows from Serre's classical theorem that $\pi_{k}(S^2)$
is finite for $k > 3$ and Theorems \ref{thm:two.sphere} and
\ref{thm:two.sphere.returned} proven in \cite{BCWW}.

Do natural features of the braid groups imply this result?
\end{enumerate}

\item Let $F_n$ denote the image of $\Theta_n(F_n)$.
Observe that the groups $\QBrun_{n+2}\cap F_{n+1}$, and
$\Brun_{n+1}\cap F_{n}$ are free. Furthermore, there is a short
exact sequence of groups
$$ 1 \to F_{n}\cap d_0(\QBrun_{n+2}) \to F_{n}\cap \Brun_{n+1} \to \pi_{n+1}S^2\to 1$$
as well as isomorphisms $$F_{n}\cap \Brun_{n+1}/(F_{n}\cap
d_0(\QBrun_{n+2})) \to \pi_{n+1}S^2$$ by 
\fullref{thm:homotopy.and.embedding braids}.

Consider the Serre 5--term exact sequence for the group extension
directly above to obtain information about the induced surjection
$$H_1(F_{n}\cap \Brun_{n+1}) \to \pi_{n+1}(S^2).$$ This
5--term exact sequence specializes to
$$
\begin{array}{ccccc}
H_2(\pi_{n+1}(S^2))&\to&H_1(F_n\cap
d_0(\QBrun_{n+2}))_{\pi_{n+1}(S^2)} &&\\
&\to &H_1(F_n\cap
\Brun_{n+1})& \to& \pi_{n+1}(S^2)\\
\end{array}
$$
where $A_{\pi}$ denotes the group of co-invariants of a $\pi$--module
$A$. Thus $\pi_{n+1}(S^2)$ is a quotient of the free abelian group
$H_1(F_n\cap \Brun_{n+1})$ with relations given by the image of the
coinvariants $H_1(F_n)\cap d_0(\QBrun_{n+2}))_{\pi_{n+1}(S^2)}$.

Give combinatorial descriptions of the induced map on the level of
the first homology groups $$H_1(F_n\cap d_0(\QBrun_{n+2})) \to
H_1(F_n\cap \Brun_{n+1}).$$
A similar problem arises with the the epimorphism $\Brun_{n+1}( S^2)
\to \pi_{n}S^2$ with kernel in the image of $\Brun_{k+1}(\mathbb
R^2)$ for $n+1 \geq 5$.

\item This problem addresses a similarity between Vassiliev
invariants and modular forms in the sense that both can be regarded
as functions defined on certain braid groups. The initial ingredient
here is given by the natural epimorphism $$B_3 \to SL(2,\mathbb Z)$$
with kernel given by the integers $\Z$. A generator for this kernel
is equal to twice a generator of the center of $B_3$.

Recall that Shimura \cite{Shimura} gives isomorphisms
\begin{gather*}H^1(SL(2,\mathbb Z);\mathbb R[x_1,x_2]) \to \mathcal
M,\\\tag*{\hbox{\hspace{28pt}and}}H^1(SL(2,\mathbb Z);\Sym^k(x_1,x_{2})) \to \mathcal
M_{2k+2}
\end{gather*}with the following properties.
\begin{enumerate}
\item The action of $SL(2,\mathbb
Z)$ on $\mathbb R[x_1,x_2]$ is specified by the tautological
representation on the two dimensional vector space with basis
$x_1,x_2$ extended multiplicatively to the polynomial ring $\mathbb
R[x_1,x_2]$.
  \item The module $\Sym^k(x_1,x_{2})$ denotes the $SL(2,\mathbb Z)$--submodule
  of $\mathbb R[x_1,x_2]$ given by polynomials of classical degree
  $k$.
\item The vector space $\mathcal M_{2k+2}$ denotes the summand of the ring of
  modular forms of weight $2k+2$ (with Shimura's weight convention).
\end{enumerate}

Also recall that the first cohomology group $H^1(SL(2,\mathbb
Z);\mathbb R[x_1,x_2])$ is given by the quotient of the module of
crossed homomorphisms $SL(2,\mathbb Z)\to \mathbb R[x_1,x_2]$ modulo
principal crossed homomorphisms. Thus elements in the classical ring
of modular forms can be regarded as equivalence classes of certain
functions defined on $SL(2,\mathbb Z)$.

Furthermore, there are isomorphisms $$H^1(B_3;
\Sym^k(x_1,x_{2})) \to E[u] \otimes \mathcal M_{2k+2}$$ where $E[u]$
denotes an exterior algebra with $u$ of degree 1. It is natural to
ask whether and how crossed homomorphisms representing
$H^1(B_3;\Sym^k(x_1,x_{2}))$ distinguish braids in $B_3$.

There is also a natural analogue for subgroups of $B_{2g+2}$ as
follows. Determine $H^1(B_{2g+2};\mathbb R[x_1, \cdots,x_{2g}])$
where $B_{2g+2}$ acts via a natural symplectic representation on a
vector space with basis $\{x_1, \cdots,x_{2g}\}$ the generating
module for the polynomial ring and $\mathbb R[x_1, \cdots,x_{2g}]$.
How do the natural crossed homomorphisms distinguish braids?

\item Consider Brunnian braids $\Brun_{k}$. Fix a braid $\gamma$
with image in the $k$-th symmetric group $\Sigma_k$ given by a
$k$--cycle. For any braid $\alpha$ in $\Brun_{k}$, the braid closure
of $\alpha \circ \gamma$ is a knot. Describe features of these knots
or those obtained from the analogous constructions for
$\Theta_k(F_{k-1})\cap \Brun_{k}$. Where do these fit in Budney's
description of the space of long knots \cite{Budney}?

\item Let $L[V]$ denote the free Lie algebra over the integers
generated by the free abelian group $V$. Let $\Der(L[V])$ denote the
Lie algebra of derivations of $L[V]$ and consider the classical
adjoint representation
$$\mathrm{Ad}\co L[V] \to \Der(L[V]).$$
Recall that the map $\Theta_k\co F_k \to P_{k+1}$ induces a
monomorphism of Lie algebras $\gr_*(\Theta_k) \co \gr_*(F_k) \to
\gr_*(P_{k+1})$ where $\gr_*(F_k)$ is isomorphic to the free Lie
algebra $L[V_k]$ with $V_k$ a free abelian group of rank $k$. In
addition, properties of the "infinitesimal braid relations" give a
representation $$\rho_k\co \gr_*(P_{k+1})\to \Der(L[V_k])$$
appearing in work on certain Galois groups \cite{Ihara} (with the
integers $\Z$ replaced by the profinite completion of $\Z$) and,
when restricted to the integers $\Z$, also addressed in
\cite{cohen.prassidis}.

Identify $F_k$ with $\Theta_k(F_k)$ in what follows below.

Give methods to describe combinatorial properties of the composite
\[
\begin{CD}
\gr_*(F_k\cap \Brun_{k+1})  @>{\gr_*(i_k)}>> \gr_*(F_k)
@>{\gr_*(\Theta_k)}>> \gr_*(P_{k+1}) @>{\rho_k}>> \Der(L[V_k])
\end{CD}
\] where $i_k\co F_k\cap \Brun_{k+1} \to F_k$ is the natural inclusion.
Let $\Phi_{k+1}$ denote this composite. When restricted to
$\gr_*(F_k)= L[V_k]$, this map is a monomorphism. Give methods to
describe the sub-quotient $$\Phi_{k+1}(\gr_*(F_k\cap
\Brun_{k+1}))/\Phi_{k+1}\gr_*(F_{k}\cap d_0(\QBrun_{k+2})).$$

\item Assume that the pure braid groups $P_n(S)$ are replaced by either
their pro-finite completion $\widehat{P_n(S)}$ or their pro-$\ell$
completions. Describe the associated changes for the homotopy groups
arising in Theorems \ref{thm:two.sphere},
\ref{thm:two.sphere.returned}, or \ref{thm:Embedding braids}. For
example, is the torsion in these homotopy groups left unchanged by
replacing $P_n(S)$ by $\widehat{P_n(S)}$?

\end{enumerate}

\section{Acknowledgements} Some of the work on this project was accomplished at the
University of Tokyo during 2004. The first author was partially
supported by the National Science Foundation under Grant No. 9704410 as
well as Darpa  Grant No. 2006-06918-01, and the Institute for
Advanced Study. The second author was partially supported by the
Academic Research Fund of the National University of Singapore No.
R-146-000-048-112.

\bibliographystyle{gtart}
\bibliography{link}

\end{document}